\documentclass[11pt, reqno]{amsart}

\usepackage{amsmath, amsthm}

\newcommand{\R}{\mathbb{R}}
\newcommand{\norm}[1]{\|#1\|}
\newcommand{\abs}[1]{|#1|}
\newcommand{\card}[1]{|#1|}
\newcommand{\epsi}{\varepsilon}

\theoremstyle{plain}
\newtheorem{theorem}{Theorem}
\newtheorem{lemma}{Lemma}

\newtheorem{proposition}{Proposition}

\begin{document}

\title{A problem of Kusner on equilateral sets}
\author{Konrad J. Swanepoel}
\thanks{This material is based upon work supported by the 
National Research Foundation under Grant number 2053752.}
\address{Department of Mathematics, Applied Mathematics and Astronomy\\ 
        University of South Africa \\ PO Box 392 \\ 
        Pretoria 0003 \\ South Africa}
\email{swanekj@unisa.ac.za}
\date{\today}

\begin{abstract}
R.~B.~Kusner [R.~Guy, Amer.\ Math.\ Monthly \textbf{90} (1983), 196--199] asked whether a set of vectors in $\R^d$ such that the $\ell_p$ distance between any pair is $1$, has cardinality at most $d+1$.
We show that this is true for $p=4$ and any $d\geq 1$, and false for all $1<p<2$ with $d$ sufficiently large, depending on $p$.
More generally we show that the maximum cardinality is at most $(2\lceil p/4\rceil-1)d+1$ if $p$ is an even integer, and at least $(1+\epsi_p)d$ if $1<p<2$, where $\epsi_p>0$ depends on $p$.
\end{abstract}

\maketitle

\section{Introduction}
Let $1<p<\infty$ and $d\geq 1$.
By $\ell_p^d$ we denote $\R^d$ endowed with the \emph{$\ell_p$-norm}
\[ \norm{x}_p=\norm{(x_1,\dots,x_d)}_p=(\sum_{i=1}^d\abs{x_i}^p)^{1/p}. \]
The \emph{unit ball} and \emph{unit sphere} (\emph{unit circle} if $d=2$) of $\ell_p^d$ are the sets $\{x: \norm{x}_p\leq 1\}$ and $\{x: \norm{x}_p= 1\}$, respectively.
Note that we do not consider the cases $p=1, \infty$ in this paper.
A set $S\subset\ell_p^d$ is \emph{$\lambda$-equilateral} ($\lambda>0$) if $\norm{x-y}_p=\lambda$ for all distinct $x,y\in S$, and \emph{equilateral} if $S$ is $\lambda$-equilateral for some $\lambda>0$.
The maximum number of elements in an equilateral set in $\ell_p^d$ is denoted by $e(\ell_p^d)$.
It is well-known that $e(\ell_2^d)=d+1$.
The standard basis vectors of $\R^d$ together with some multiple of $(1,1,\dots,1)$ demonstrates that $e(\ell_p^d)\geq d+1$ for all $1<p<\infty$.
A result of Petty \cite{MR43:1051} gives as a special case that $e(\ell_p^d)<2^d$ for $d\geq 2$.
It is also well-known that $e(\ell_p^2)=3$ (see e.g\ \cite[Section~5]{MR2002h:46015a}).
Kusner \cite{Guy} asked whether $e(\ell_p^d)=d+1$ for all $d\geq 2$ and $1<p<\infty$.
This problem has recently been studied by Smyth \cite{Smyth} and Alon and Pudl\'ak \cite{AP}.
Smyth showed $e(\ell_p^d)< c p d^{(p+1)/(p-1)}$ for some $c>0$, and also $e(\ell_p^d)=d+1$ for $2-\alpha_d<p<2+\alpha_d$ where $\alpha_d=O(1/(d\log d))$ 
(this second statement also follows from a more general result of Brass \cite{MR2000i:52012} and Dekster \cite{2001b:52001}).
The general upper bound was improved by Alon and Pudl\'ak to $e(\ell_p^d)< c_p d^{(2p+2)/(2p-1)}$ for some $c_p>0$ depending on $p$.
For $p$ an even integer, Galvin (see \cite{Smyth}) showed $e(\ell_p^d)\leq 1+(p-1)d$, while for $p$ an odd integer Alon and Pudl\'ak showed $e(\ell_p^d)\leq c_p d\log d$ for some $c_p>0$.

First of all we improve Galvin's result as follows.

\begin{theorem}\label{th1}
For $p$ an even integer and $d\geq 1$ we have
\[ e(\ell_p^d) \leq \begin{cases} (\tfrac{p}{2}-1)d+1 & \text{if } p\equiv 0 \pmod{4},\\
\tfrac{p}{2}d+1 & \text{if } p\equiv 2 \pmod{4}. \end{cases} \]
In particular, $e(\ell_4^d)=d+1$.
\end{theorem}
By a compactness argument we thus have that $e(\ell_p^d)=d+1$ for $p$ in a small interval around $4$, the size of the interval depending on $d$.
We have no information on the size of this interval, nor on whether $e(\ell_p^d)=d+1$ for any other values of $p>2$.
The proof of Theorem~\ref{th1} in Section~\ref{proof1} uses a linear algebra method (see \cite[Part~III]{Jukna} for an exposition).

Secondly we show that $e(\ell_p^d)>d+1$ holds for all $1<p<2$ if $d$ is sufficiently large.

\begin{theorem}\label{th2}
For any $1<p<2$ and $d\geq 1$, let
\[k=\left\lceil\frac{\log(1-2^{p-2})^{-1}}{\log 2}\right\rceil-1.\]
Then
\[e(\ell_p^d)\geq \left\lfloor\frac{2^{k+1}}{2^{k+1}-1}d\right\rfloor.\]
In particular, if $d\geq 2^{k+2}-2$ then $e(\ell_p^d)>d+1$.
\end{theorem}
For example, if $1<p\leq\frac{\log 3}{\log 2}$, then $e(\ell_p^d)\geq\lfloor 4d/3\rfloor$ and $e(\ell_p^6)\geq 8$.
For $p$ close to $2$ the theorem implies that
$e(\ell_p^d)>d+1$ if $p<2-\Omega(1/d)$.
Thus we have reached the above-mentioned bound of Smyth except for a $\log d$ factor.
Theorem~\ref{th2} is proved in Section~\ref{proof2} by constructing explicit examples based on Hadamard matrices.

The smallest dimension for which Theorem~\ref{th2} gives an example of $e(\ell_p^d)>d+1$ is $d=6$.
With a slightly modified construction we also give examples for $d=4$.
However, we have no examples for $d=3$ or $d=5$.
\begin{theorem}\label{th3}
For any $1<p\leq\frac{\log5/2}{\log 2}$ we have $e(\ell_p^4)\geq 6$.
\end{theorem}
The proof is also in Section~\ref{proof2}.

\section{Upper bounds for $p$ an even integer}\label{proof1}
Let $S$ be a $1$-equilateral set in $\ell_p^d$ where $p$ is an even integer.
For each $a\in S$, let $P_a(x)=P_a(x_1,\dots,x_d)$ be the following polynomial:
\begin{eqnarray}
 P_a(x) &:=& -1 +\norm{x-a}_p^p \notag\\
&=& -1+\norm{a}_p^p + \sum_{i=1}^d x_i^p +\sum_{i=1}^d\sum_{m=1}^{p-1} \binom{p}{m}(-a_i)^{p-m}x_i^m. \label{eq1}
\end{eqnarray}
Thus each $P_a$ is in the linear span of
\[ \{1,\sum_{i=1}^d x_i^p\}\cup\{x_i^m: 1\leq m\leq p-1;\, 1\leq i\leq d\}, \]
which is a subspace of dimension $(p-1)d+2$ of the vector space of real polynomials in the variables $x_1,\dots,x_d$.
Since $P_a(a)=-1$ for all $a\in S$ and $P_a(b)=0$ for all distinct $a,b\in S$, we have that $\{P_a:a\in S\}$ is linearly independent.
Thus we already have $\card{S}\leq (p-1)d+2$.
We now show that the larger set
\[ \mathcal{P} := \{P_a:a\in S\} \cup\{1\}\cup\{x_i^m: 1\leq i\leq d;\, 1\leq m\leq k\} \]
is still linearly independent, where $k=p/2$ if $p\equiv 0 \pmod{4}$ and $k=p/2-1$ otherwise.
This will give $\card{S}+1+kd\leq (p-1)d+2$, proving Theorem~\ref{th1}.

We only consider the case $p\equiv0\pmod{4}$, the other case being similar.
Let $\lambda, \lambda_a\, (a\in S), \lambda_{i,m} \, (1\leq i\leq d;\, 1\leq m\leq p/2)$ be real numbers satisfying
\begin{equation} 
\lambda 1 + \sum_{a\in S}\lambda_a P_a(x) + \sum_{i=1}^d\sum_{m=1}^{p/2}\lambda_{i,m}x_i^m \equiv 0.\label{eq2}
\end{equation}
If we substitute \eqref{eq1} into \eqref{eq2} we obtain
\begin{multline}
 \lambda 1 + \sum_{a\in S}\lambda_a(-1+\norm{a}_p^p) + \sum_{i=1}^d(\sum_{a\in S}\lambda_a)x_i^p \\
+ \sum_{i=1}^d\sum_{m=1}^{p-1}\sum_{a\in S}\lambda_a\binom{p}{m}(-a_i)^{p-m}x_i^m + \sum_{i=1}^d\sum_{m=1}^{p/2}\lambda_{i,m}x_i^m \equiv 0.
\end{multline}
Thus the coefficients of this polynomial are all $0$, giving
\begin{equation}\label{eq3}
\lambda+\sum_{a\in S}\lambda_a(-1+\norm{a}_p^p)=0,
\end{equation}
\begin{equation}\label{eq4}
\sum_{a\in S}\lambda_a=0,
\end{equation}
\begin{equation}\label{eq6}
\lambda_{i,m}+\sum_{a\in S}\lambda_a\binom{p}{m}(-a_i)^{p-m}=0 \qquad \forall m=1,\dots,p/2;\, i=1,\dots, d,
\end{equation}
\begin{equation}\label{eq5}
\sum_{a\in S}\lambda_a a_i^m = 0 \qquad \forall m=1,\dots,p/2-1.
\end{equation}
Substitute $x=b\in S$ into \eqref{eq2}:
\begin{equation}\label{eq7}
-\lambda_b + \lambda + \sum_{i=1}^d\sum_{m=1}^{p/2}\lambda_{i,m}b_i^m=0\quad\forall b\in S.
\end{equation}
Multiply \eqref{eq7} by $-\lambda_b$ and sum over all $b\in S$:
\[ \sum_{b\in S}\lambda_b^2 - \lambda\sum_{b\in S}\lambda_b - \sum_{i=1}^d\sum_{m=1}^{p/2}\lambda_{i,m}\sum_{b\in S}\lambda_b b_i^m = 0.\]
By \eqref{eq4} and \eqref{eq5} this simplifies to
\[ \sum_{b\in S}\lambda_b^2 - \sum_{i=1}^d \lambda_{i,p/2}\sum_{b\in S}\lambda_b b_i^{p/2} = 0,\]
which by \eqref{eq6} simplifies to
\[ \sum_{b\in S}\lambda_b^2 + \sum_{i=1}^d\binom{p}{m}(\sum_{a\in S}\lambda_a a_i^{p/2})^2 = 0.\]
Since the left-hand side is a sum of squares, $\lambda_b=0$ for all $b\in S$. 
By \eqref{eq3} we then have $\lambda=0$, and by \eqref{eq6} $\lambda_{i,m}=0$ for all $m$ and $i$.
Thus the set $\mathcal{P}$ is linearly independent, finishing the proof.
\qed

\section{Lower bounds for $1<p<2$}\label{proof2}
According to the following proposition, if we can find a $2^{1/p}$-equilateral set of $k+1$ points on the unit sphere of $\ell_p^k$, we can construct equilateral sets in $\ell_p^d$ of more than $d+1$ points if $d$ is sufficiently large.
The construction is similar to the Lenz construction in combinatorial geometry (see \cite[pp.~148, 159, 194]{MR96j:52001}).
\begin{proposition}\label{prop1}
Let $1<p<\infty$, $k, d\geq 1$.
If $\ell_p^k$ has a $2^{1/p}$-equilateral set of cardinality $k+1$ on the unit sphere, then $\ell_p^d$ has a $2^{1/p}$-equilateral set of cardinality $\lfloor(1+\frac{1}{k})d\rfloor$ on the unit sphere.
\end{proposition}
\vspace{-5mm}
\begin{proof}
Let $m=\lfloor d/k \rfloor$ and $r=d-km$.
Let $\ell_p^d=\overbrace{\ell_p^k\oplus\dots\oplus\ell_p^k}^{m \text{ times}}\oplus\ell_p^r$.
Let the equilateral set in $\ell_p^k$ be $S=\{v_1,\dots,v_{k+1}\}$.
Let $S_i$ be the copy of $S$ in the $i$'th copy of $\ell_p^k$ in $\ell_p^d$, $i=1,\dots,m$, and let $S_0$ be the copy of the standard unit vectors $e_1,\dots,e_r$ in the copy of $\ell_p^r$, which is also a $2^{1/p}$-equilateral set of unit vectors.
Clearly the distance between a vector in $S_i$ and a vector in $S_j$ is $2^{1/p}$ for distinct $i,j$, since both are unit vectors.
Thus $S_0\cup\dots\cup S_m$ is the required set, since it has cardinality $m(k+1)+r=d+m=d+\lfloor d/k\rfloor$.
\end{proof}

Before we construct the required $2^{1/p}$-equilateral sets, we need a technical two-dimensional result.
\begin{lemma}\label{lemma1}
Let $1<p<2$.
For each $\lambda\in [2^{1-1/p}, 2^{1/p}]$ there exist unit vectors $u, v\in\ell_p^2$ such that $\norm{u+v}_p=\norm{u-v}_p=\lambda$.
\end{lemma}
Geometrically the lemma says that there exists a quadrilateral inscribed in the unit circle of $\ell_p^2$ with all four sides of length $\lambda$, for any $\lambda\in [2^{1-1/p}, 2^{1/p}]$.
This is easily seen for $\lambda=2^{1/p}$ ($u=(1,0)$ and $v=(0,1)$) and for $\lambda=2^{1-1/p}$ ($u=((1/2)^{1/p},(1/2)^{1/p})$ and $v=(-(1/2)^{1/p},(1/2)^{1/p})$.
The inbetween values are then covered by a continuity argument.
We omit the details.

We also need Hadamard matrices.
Recall that a Hadamard matrix of order $k$ is a $k\times k$ matrix $H$ with all entries $\pm 1$, satisfying $HH^T=kI$.
It is well-known that if a Hadamard matrix of order $k$ exists, then $k=1,2$ or $k$ is divisible by $4$ \cite{MR2002i:05001}, and not known whether the converse holds.
However, we only need the fact that Hadamard matrices of order $2^n$ exist for all $n\geq 0$, as shown by the well-known inductive construction
\[ H_0 = [1], \qquad H_{n+1}=\left[\begin{array}{cc} H_n & H_n \\ H_n & - H_n\end{array}\right], n\geq 0.\]

\begin{proposition}\label{prop2}
Suppose that there exists a Hadamard matrix of order $k\geq 2$.
Then for any
\[ p\in \left[2+\frac{\log(1-k^{-1})}{\log 2}, 2+\frac{\log(1-(2k)^{-1})}{\log 2}\right], \quad p\neq 1,\]
there exists a $2^{1/p}$-equilateral set of cardinality $2k$ on the unit sphere of $\ell_p^{2k-1}$.
\end{proposition}
\begin{proof}
We may assume without loss of generality that the $k\times k$ Hadamard matrix is normalized such that its first column contains only $+1$'s.
Delete the first column and let the $k$ resulting rows be $w_1,\dots, w_k\in\{\pm1\}^{k-1}$.
We have that any two distinct $w_i$ and $w_j$ differ in exactly $k/2$ coordinates.
Let
\begin{equation}\label{eqq1}
\lambda=2\left(\frac{(3-2^{p-1})k-2}{2(k-1)}\right)^{1/p}.
\end{equation}
The given bounds on $p$ ensure that $2^{1-1/p}\leq \lambda\leq 2^{1/p}$.
We now take $u,v\in\ell_p^2$ from Lemma~\ref{lemma1}, and let
\[ u_i=(\mu,w_i\otimes u),\, v_i=(-\mu,w_i\otimes v)\in\ell_p^{2k-1}, \quad i=1,\dots,k,\]
where
\begin{equation}\label{eqq2}
\mu=((2^{p-2}-1)k+1)^{1/p}
\end{equation}
and the \emph{Kronecker product} $a\otimes b$ for any vectors $a=(a_1,\dots,a_m)\in\R^m$ and $b\in\R^n$ is defined as $(a_1 b, a_2 b, \dots, a_m b)\in\R^{mn}$.
The given lower bound on $p$ ensures that $\mu$ is well-defined.
Then for any distinct $i, j$ we have $\norm{u_i-u_j}_p^p=(k/2)\norm{2u}_p^p=2^{p-1}k$ and similarly $\norm{v_i-v_j}_p^p=2^{p-1}k$.
Also for any $i,j$,
\begin{eqnarray*}
\norm{u_i-v_j}_p^p &=& (2\mu)^p + \sum_{m=1}^{k-1}\norm{u+\epsi_m v}_p^p\quad \text{for some } (\epsi_m)\in\{\pm 1\}^{k-1}\\
&=& (2\mu)^p + (k-1)\lambda^p \quad \text{by Lemma~\ref{lemma1}}\\
&=& 2^{p-1}k \quad \text{by \eqref{eqq1} and \eqref{eqq2}.} 
\end{eqnarray*}
Thus $S:=\{u_1,\dots,u_k, v_1,\dots, v_k\}$ is equilateral.
Also,
\[\norm{u_i}_p^p=\mu^p+(k-1)\norm{x}_p^p=\mu^p+k-1=2^{p-2}k,\]
by \eqref{eqq2}, and similarly, $\norm{v_i}_p^p=2^{p-2}k$.
Thus if we scale $S$ by $(2^{p-2}k)^{-1/p}$, we obtain a $2^{1/p}$-equilateral set of unit vectors of cardinality $2k$.
\end{proof}

Note that the two smallest dimensions $d$ for which the above proposition ensures a $2^{1/p}$-equilateral set of unit vectors of size $d+1$ are $d=3$ (with $1<p\leq\log3/\log 2$) and $d=7$ (with $\log 3/\log 2\leq p\leq \log(7/2)/\log 2$).
It is not difficult to see that such a set does not exist for $d=2$ and any $1<p<\infty$, and also not for $p=2$ and any $d$.
We do not know whether such sets exist if $d=4, 5, 6$.
It is doubtful that they exist for $p>2$.

\begin{proof}[Proof of Theorem~\ref{th2}]
Note that the given value of $k$ ensures that 
\[ 2+\frac{\log(1-2^{-k})}{\log 2} < p \leq 2+\frac{\log(1-2^{-k-1})}{\log 2},\]
and thus Proposition~\ref{prop2} applied to a Hadamard matrix of order $2^k$ gives a $2^{1/p}$-equilateral set of cardinality $2^{k+1}$ on the unit sphere of $\ell_p^{2^{k+1}-1}$.
Then by Proposition~\ref{prop1} $\ell_p^d$ has a $2^{1/p}$-equilateral set of size $\lfloor (1+(2^{k+1}-1)^{-1})d\rfloor$.
\end{proof}

\begin{proof}[Proof of Theorem~\ref{th3}]
Let $\lambda=2(3-2^p)^{1/p}$ and take $u,v\in\ell_p^2$ from Lemma~\ref{lemma1}.
(Since $1<p\leq\log(5/2)/\log 2$ we have $2^{1-1/p}\leq \lambda \leq 2^{1/p}$.)
Let $\mu=(2^p-2)^{1/p}$.
Then a simple calculation shows that
\[S=\{(\mu,\pm u, 0),(-\mu,\pm v, 0),(0,o,\pm 1)\}\subset\ell_p^4\]
 is equilateral.
\end{proof}

\providecommand{\MR}{\relax\ifhmode\unskip\space\fi MR }
\providecommand{\MRhref}[2]{%
  \href{http://www.ams.org/mathscinet-getitem?mr=#1}{#2}
}
\providecommand{\href}[2]{#2}


\begin{thebibliography}{10}

\bibitem{AP}
N.~Alon and P.~Pudl{\'a}k, \emph{Equilateral sets in $l_p^n$}, to appear in
  Geometric and Functional Analysis.

\bibitem{MR2000i:52012}
P.~Bra{\ss}, \emph{On equilateral simplices in normed spaces}, Beitr\"age
  Algebra Geom.\ \textbf{40} (1999), 303--307. \MR{2000i:52012}

\bibitem{2001b:52001}
B.~V.~Dekster, \emph{Simplexes with prescribed edge lengths in Minkowski and 
Banach spaces}, Acta Math.\ Hungar.\ \textbf{86} (2000), 343--358. \MR{2001b:52001}


\bibitem{Guy}
R.~K.~Guy, \emph{An olla-podrida of open problems, often oddly posed}, Amer.\
  Math.\ Monthly \textbf{90} (1983), 196--199.

\bibitem{Jukna}
S.~Jukna, \emph{Extremal combinatorics}, Springer-Verlag, Berlin-Heidelberg-New
  York, 2001.

\bibitem{MR2002h:46015a}
H.~Martini, K.~J.~Swanepoel, and G.~Wei{\ss}, \emph{The geometry of
  {M}inkowski spaces---a survey. {I}}, Expo.\ Math.\ \textbf{19} (2001),
  97--142. \MR{2002h:46015a}. Errata: Expo.\ Math.\ \textbf{19} (2001), p.~364. 
  \MR{2002h:46015b}

\bibitem{MR96j:52001}
J.~Pach and P.~K. Agarwal, \emph{Combinatorial geometry},
  John Wiley \& Sons Inc., New York, 1995.
  \MR{96j:52001}

\bibitem{MR43:1051}
C.~M. Petty, \emph{Equilateral sets in {M}inkowski spaces}, Proc.\ Amer.\ Math.\
  Soc.\ \textbf{29} (1971), 369--374. \MR{43 \#1051}

\bibitem{Smyth}
C.~Smyth, \emph{Equilateral or $1$-distance sets and {Kusner's} conjecture},
  submitted.

\bibitem{MR2002i:05001}
J.~H. van Lint and R.~M. Wilson, \emph{A course in combinatorics}, second ed.,
  Cambridge University Press, Cambridge, 2001. \MR{2002i:05001}

\end{thebibliography}
\end{document}